\documentclass[11pt]{article}
\usepackage{amssymb,latexsym,amsmath,amsbsy,amsthm,amsxtra,amsgen,dsfont,pdfsync,color}

\newtheorem{thm}{Theorem}
 
\newtheorem{prop}{Proposition}

\begin{document}
\title{Fubini's Theorem for Daniell Integrals}

\author{G\"otz Kersting\thanks{Institut f\"ur Mathematik, Goethe Universit\"at, Frankfurt am Main, Germany,  kersting@math.uni-frankfurt.de} ,  Gerhard Rompf}

\maketitle

\begin{abstract}
We show that in the theory of Daniell integration iterated integrals may always be formed, and the order of integration may  always be interchanged. By this means, we discuss product integrals and show that the related Fubini theorem holds in full generality. The results build on a density theorem on Riesz tensor products due to Fremlin, and on the Fubini-Stone Theorem.

\vspace{.2cm}

\vspace{.2cm}

\end{abstract}

\section{Introduction and main results}

The Daniell  integral is considered to be a particulary elegant approach to integration. Dating from 1918 \cite{Dan}, it took, nonetheless, several decades, until multiple Daniell integration had been treated in the literature. Two publications stand out. Stone \cite{Sto2} showed in 1948 for an integral on a product space, which can be reduced to an iterated integral,   that its extension to summable functions may as well be expressed by an iterated integral. By contrast, Fremlin's far reaching results from 1972 \cite{Fre} allow  the  construction of  a product of integrals without touching iterated integrals. In this note we like to meet both these aspects and to proceed to a Fubini type theorem in such generality, as it is long-established in the narrower domain of Lebesgue integration. Our main observation is that  iterated integrals may  be formed for any two Daniell integrals, and that the order of integration may always be interchanged.

Let us recall  some major notions from Daniell integration theory. Let $F(Z)$ be the set of all function from some non-empty set~$Z$ into the real numbers. Consider   a set $L\subset F(Z)$, $L\neq \emptyset$, and a mapping $I:L\to \mathbb R$. If the conditions

\begin{enumerate}
\item[(i)] (Additivity) If $f,g\in L$, then $f+g\in L$ and $I(f+g)=I(f)+I(g)$.
\item[(ii)] (Homogeneity) If $f \in L$ and $r \in \mathbb R$, then $rf\in L$ and $I(rf)=rI(f)$.
\item[(iii)] (Monotonicity) If $f\in L$, then $|f|\in L$ and $I(f)\le I(|f|)$.
\item[(iv)] (Continuity) If $\inf_n f_n =0$ for a decreasing sequence $f_1 \ge f_2\ge \cdots$  in $L$, then $\inf_n I(f_n)=0$.
\end{enumerate}
are fulfilled, then one calls $L$ a {\em Riesz space} (vector lattice) over $Z$ and $I$ an  {\em elementary Daniell integral} on $L$.   Then, together with $f$ and $g$ also $\sup(f,g)$ and $\inf(f,g)$ belong to $L$. (Here all operations or relations of functions are pointwisely defined.) The smallest Riesz space within $F(Z)$ containing some set of functions $S \subset F(Z)$ is denoted by $L(S)$. 
We take the liberty to call an elementary Daniell integral $I$ briefly an {\em integral}. 
We recall from the celebrated integration theory of Daniell,  as presented e.g. in \cite{Bla,Loo,Roy,Tay} that any integral $I:L\to \mathbb R$ may be extended to the  integral $ I: L^1 \to \mathbb R$, where $L^1\supset L$ is the Riesz space of all real-valued $I$-summable  functions. To this extension we may apply the theorems of monotone  and dominated convergence, as known from Lebesgue integration (see e.g. \cite{Roy,Tay}). (As customary, we denote the  integral and its extension by one and the same symbol.)

Now, let us consider two integrals $J:G\to \mathbb R$ and $K:H\to \mathbb R$ with Riesz spaces $G$ and $H$ over $X$ and $Y$, respectively. As indicated above, 
in order to combine them there are two aspects. First, by analogy to Lebesgue product measures, one  asks for a product integral $I:L\to \mathbb R$ over  $X\times Y$ fulfilling $I(g\otimes h)=J(g)K(h)$ for any $g\in G^1$, $h \in H^1$, with $g\otimes h(x,y):=g(x)h(y)$ for all $(x,y)\in X\times Y$. The natural candidate for $L$ is the space $L(G^1\otimes H^1)$, the smallest Riesz space over $X\times Y$  containing the tensor product $G^1\otimes H^1$ (the vector space spanned by all $g\otimes h$ with $g\in G^1$ and $h\in H^1$). This Riesz space $L$ serves as the Riesz tensor product of $G^1$ and $H^1$ (characterized by an universality property, see \cite{Fre}).

Second, one would like to express product integrals by iterated integrals. Let  us consider  for some  function $f \in F(X\times Y)$ and the derived functions $f^x:Y\to \mathbb R$,  $f^y:X\to \mathbb R$,  given by $f^x(y)= f^y(x)=f(x,y)$,  the conditions

\begin{enumerate}
\item[(a)] $f^x\in H^1$  and $f^y \in G^1$ for all $x \in X$, $y\in Y$,
\item[(b)] $K(f)\in G^1$ and $J(f)\in H^1$, with $K(f):X \to \mathbb R$ and $J(f):Y\to \mathbb R$ given by $K(f)(x):=K(f^x)$, $J(f)(y):=J(f^y)$.
\end{enumerate}
Then the iterated integrals $J(K(f))$ and $K(J(f))$ are well-defined.

\begin{thm} \label{thm1} For any two integrals $J:G\to \mathbb R$ and $K:H\to \mathbb R$ there is a unique integral $I:L(G^1\otimes H^1) \to \mathbb R$ fulfilling for all $g\in G^1$, $h\in H^1$ the equation 
\[ I(g\otimes h)=J(g)\cdot K(h) .\]
Moreover,  any $f\in L(G^1\otimes H^1)$ meets the conditions {\em (a)} and {\em (b)}, and we have
\[ I(f)= J(K(f))=K(J(f)).\]
\end{thm}

\noindent 
 This integral $I$ is written $I=J\otimes K$ and called the {\em product integral of $J$ and} $K$.

With this theorem at hand, it is straightforward to proceed  to a related Fubini type result. To this end, we somewhat weaken the requirements (a) and (b) and consider for $f\in F(X\times Y)$  the conditions
\begin{enumerate}
\item[(a')]  there are a $J$-null set  $A\subset X$ and a $K$-null set $B\subset Y$ such that $f^x\in H^1$ and $f^y \in G^1$  for $x\notin A$ and $y\notin B$,
\item[(b')] there are functions $g\in G^1$ and $h\in H^1$ such that $g(x)=K(f^x)$ and $h(y)=J(f^y)$ for  $x\notin A$ and  $y\notin B$.
\end{enumerate}
We recall that, for an integral $I$ on the Riesz space $L$ over $Z$, a set $A\subset Z$ is called an $I$-null set, if  its characteristic function $1_A$ has the property: If $f\in L$, then $f\cdot 1_A\in L$ and $I(f\cdot 1_A)=0$. Note that the above functions $g$ and $h$ are in general not uniquely determined, yet the values of  $J(g)$ and $K(h)$ are unique. Thus, it is justified to  write $K(f)$ and $J(f)$ instead of $g$ and $h$, respectively, and to build  the iterated integrals $J(K(f))$ and $K(J(f))$ for  functions $f$ satisfying (a') and (b').

Now our Fubini type theorem reads as follows.

\begin{thm} \label{thm2} Let $I: L \to \mathbb R$ be the product integral of $J:G\to \mathbb R$ and $K:H\to \mathbb R$. Then any function $f \in L^1$ fulfils the conditions {\em (a')} and {\em (b')}, and we have
\[ I(f)= J(K(f))=K(J(f)).\]
\end{thm}

\noindent
Our theorems have been obtained in the literature for different special cases, as in Loomis \cite{Loo} or Zaanen \cite{Zaa}. In these cases the tensor products \mbox{$G\otimes H$} are already Riesz spaces, which allows a direct construction of the product integral. Also, in these cases the Riesz spaces under consideration  fulfil the Stone condition ($f\in L$ entails $\min(f,1)\in L$), implying that one may re\-present the Daniell integrals by Lebesgue integrals and draw on the classical Fubini theorem, too. In the general case we are mainly faced with the task to bridge the gap from $G^1\otimes H^1$ to $L(G^1\otimes H^1)$.

The proofs, given in the next section, rely on two results from the literature. Theorem \ref{thm1} hinges on a density result, concerning  the inclusion  $G^1\otimes H^1\subset L(G^1\otimes H^1)$.  Fremlin obtained this remarkable result in much greater generality \cite[Theorem 4.2 (iii)]{Fre}. His proof   is   extensive and uses advanced representation theorems for Riesz spaces. A more intrinsic derivation  was given by Grobler and Labuschagne \cite{Gro}, still their approach is somewhat involved and demanding. From these considerations only a smaller part  is needed in our case, which we present   in the subsequent selfcontained Proposition \ref{prop1}  (with a view to  readers with little experience in the theory of Riesz spaces). This density statement allows to extent the product measure $I$ from $G^1\otimes H^1$ to the enclosing Riesz space.

We note that the first half of Theorem \ref{thm1} may as well be derived with marginal efforts from \cite[Theorem 5.3]{Fre}. There, Fremlin constructs a linear functional $I$ on the Riesz tensor product $L$ of $G$ and $H$ by approximating $I(f)$, $f\in L$, from above and from below by the terms $I(u)$ with $u\in G \otimes H$ and $u\ge f$ or $u\le f$, respectively. We shall use   iterated integrals for this purpose.

Theorem \ref{thm2} is a direct consequence of the second part of Theorem \ref{thm1} together with   the Fubini-Stone theorem \cite{Sto2}, which  found its way into text books like \cite{Tay,Zaa}. Below, we recall it in Proposition \ref{prop2} (see also \cite{Ric, Mik,Dia}).

\section{Proofs}

The following proposition is a special case of Fremlin \cite[Theorem 4.2]{Fre}.  Freudenthal's spectral theorem \cite[Theorem 40.2]{LuZa} could be used here, too. Our comparatively elementary proof borrows from Grobler et al \cite{Gro} and is much in the spirit of the spectral theorem. 

Recall that a Riesz space $L$ is called {\em Dedekind $\sigma$-complete}, if it has the property
\[ \text{If } f,f_1,f_2, \ldots \in L \text{ fulfil } f_1\le f_2 \le  \cdots \le f, \text{ then } \sup_n f_n \in L. \]

\begin{prop} \label{prop1}  Let $G$ and $H$ be Dedekind $\sigma$-complete Riesz spaces over $X$ and $Y$. Then for every $f\in L(G\otimes H)$ there is a pair $(g,h)\in G^+\times  H^+$ and for every $\varepsilon>0$ some $u_\varepsilon \in G\otimes H$ such that $|f|\le g\otimes h$ and
\[ |f-u_\varepsilon| \le \varepsilon g\otimes h. \]
\end{prop}

\begin{proof}
Preliminary, we recall a well-known property of a Dedekind $\sigma$-complete Riesz space $L$ over $Z$: Let $f,u_1, \ldots, u_k,v_1, \ldots, v_k \in L$, then
\begin{align}\label{eq1} f':=f\cdot 1_{u_1 < v_1} \cdots 1_{u_k<v_k} \in L \end{align}
as well. Here $1_{u<v}$ denotes the characteristic function of the set $\{z\in Z: u(z)<v(z) \}$. For the proof it is sufficient to consider the case $f\ge 0$. Then we have $f' = \sup_n \min\{ f, n(v_1-u_1)^+, \ldots, n(v_k-u_k)^+\} \le f$, which entails equation \eqref{eq1}.

Now, for $(g,h)\in G^+\times H^+$ set $G_g:= \{v\in G: |v|\le cg \text{ for some } c>0\}$  and $H_h:=\{w\in H: |w|\le ch \text{ for some } c>0\}$. It is straightforward to see that $G_g$ and $H_h$ are Riesz spaces. Further, let 
\[L_{g,h}:= \{ f\in L  \mid \forall \varepsilon>0 \ \exists  u_\varepsilon \in G_g\otimes H_h :  |f-u_\varepsilon | \le \varepsilon g\otimes h\}.\]
The major task now is to show that $L_{g,h}$ is a Riesz space. It is easy to see that $L_{g,h}$ is a vector space, thus it remains to show that with $f\in L_{g,h}$ we also have $|f|\in L_{g,h}$.

First, we consider the case $f=\sum_{i=1}^r g_i\otimes h_i\in G_g\otimes H_h$. Obviously, $f\in L_{g,h}$. In order to show that $|f|\in L_{g,h}$, let for any $n\in \mathbb N$ 
\[ A_{a_1\ldots a_kn}:= \Big\{ \frac {a_1}n c g \le g_1 < \frac {a_1+1}nc g \Big\} \cap \cdots \cap \Big\{ \frac {a_k}n c g \le g_k < \frac {a_k+1}n c g\Big\} \]
with integers $a_1, \ldots, a_k\in\{-n, \ldots,n\}$. Apparently, for $c>0$  large enough  these sets  together with the set $\{g=0\}$ yield a partition of the domain $X$. Similar for large $c>0$ the sets
\[ B_{b_1\ldots b_kn}:= \Big\{ \frac {b_1}n c h \le h_1 < \frac {b_1+1}nc h \Big\} \cap \cdots \cap \Big\{ \frac {b_k}n c h \le h_k < \frac {b_k+1}n c h\Big\} \]
with $b_1, \ldots,  b_k \in\{-n, \ldots,n\}$ and the set $\{h=0\}$ make a partition of $Y$. From \eqref{eq1} we obtain that \mbox{$g\cdot 1_{A_{a_1\ldots a_kn}}$} belongs for all $a_1, \ldots, a_k$ to $G$, and $h\cdot 1_{B_{b_1\cdots b_kn}}$ to $H$. Define
\[ f_n := \sum_{a_1,b_1, \ldots a_k,b_k=-n}^n \Big(\frac {a_1b_1}{n^2}+ \cdots +\frac{a_kb_k}{n^2}\Big)\Big((cg\cdot 1_{A_{a_1\ldots  a_kn} })\otimes (ch\cdot 1_{B_{b_1\ldots  b_kn} })\Big).\] 
Since the sets $A_{a_1\ldots  a_kn} \times B_{b_1\ldots  b_kn} $ are pairwise disjoint, there will be for each $f_n(x,y)$ with $x\in X$ and $y\in Y$ at most one non-vanishing summand in the previous sum. Hence
\[ |f_n| =\sum_{a_1,b_1, \ldots a_k,b_k=-n}^n c^2\Big|\frac {a_1b_1}{n^2}+ \cdots +\frac{a_kb_k}{n^2}\Big|\Big((g\cdot 1_{A_{a_1\ldots  a_kn} })\otimes (h\cdot 1_{B_{b_1\ldots  b_kn} })\Big),\]
consequently $|f_n|\in G_g\otimes H_h$. Moreover, for $x\in A_{a_1\ldots  a_kn}$ and $y\in B_{b_1\ldots  b_kn} $, since $|g_i| \le cg$, 
\begin{align*}\Big| g_i(x)h_i(y)&- \frac {a_i b_i}{n^2}c^2g(x)h(y)\Big| \\&\le \Big|g_i (x)\Big(h_i(y)- \frac{b_i}nch(y)\Big)\Big|+ \Big|\frac{b_i}nch(y)\Big(g_i(x)- \frac{a_i} ncg(x)\Big) \Big|\\
&\le \frac {2c^2}n g(x)h(y). 
\end{align*}
Therefore, in view of 
\begin{align*}
f-f_n =  \sum_{a_1,b_1, \ldots a_k,b_k=-n}^n \sum_{i=1}^k\Big(g_i\otimes h_i -c^2 \frac {a_i b_i}{n^2} g\otimes h\Big) 1_{A_{a_1\ldots  a_kn} \times B_{b_1\ldots  b_kn}} ,
\end{align*}
we obtain
\[\big||f|-|f_n|\big| \le |f-f_n| \le \frac {2kc^2}n g\otimes h . \]
This estimate implies that $|f|\in L_{g,h}$.  

Now let $f$ be any element of $L_{g,h}$. Then, for any $\varepsilon >0$ there is a $u\in G_g\otimes H_h$ such that $|f-u|\le \varepsilon g\otimes h$, moreover, as just shown, there is a $v\in G_g\otimes H_h$ such that $||u|-v|\le \varepsilon g \otimes h$. Therefore
\[ ||f|-v| \le ||f|-|u|| + ||u|-v| \le |f-u|+||u|-v| \le 2\varepsilon g\otimes h . \]
Thus $|f|\in L_{g,h}$. Altogether we have shown that $L_{g,h}$ is a Riesz space.

In order to finish the proof, let $L':= \bigcup L_{g,h}$, where the union is taken over all $(g,h)\in G^+\times H^+$. Let $f_1,f_2 \in L'$, thus $f_i \in L_{g_i,h_i}$ with suitable $g_i,h_i$, $i=1,2$. Then $f_1,f_2 \in L_{g_1+g_2, h_1+h_2}$, consequently $f+g\in L'$. Also $\lambda f_1, |f_1| \in L'$, meaning that $L'$ is a Riesz space. Further, for any $(g,h)\in G\times H$ we have $g\otimes h \in L_{|g|,|h|} \subset L'$, thus $G\otimes H\subset L'$. Since $L'\subset L$ we obtain $L'=L$, which implies the proposition's claim.
\end{proof}

\newpage

\begin{proof}[Proof of  Theorem \ref{thm1}] First, we prove that $J(K(f))$ and $K(J(f))$ are well-defined for $f\in L:=L(G^1\otimes H^1)$. From dominated convergence it follows that the Riesz spaces $G^1$ and $H^1$ are Dedekind $\sigma$-complete. Thus, for any $f \in L$ there are due to Proposition \ref{prop1}
 $f_n=\sum_{i=1}^{r_n} g_{in}\otimes h_{in} \in  G^1\otimes  H^1$ such that  $|f-f_n| \le n^{-1} g\otimes h$ for some non-negative $g\in G^1$ and $h\in  H^1$. Hence, for any $x\in X$ we have 
\begin{align*} |f^x - \sum_{i=1}^{r_n} g_{in}(x)h_{in}|\le n^{-1} g(x)h. \end{align*}
Taking  the limit $n\to \infty$ and using dominated convergence yields  $f^x \in  H^1$ for all $x$, furthermore
\[ | K(f^x) - \sum_{i=1}^{r_n} g_{in}(x) K(h_{in})|\le n^{-1} g(x) K(h). \]
Again taking the limit, dominated convergence yields $ K(f) \in  G^1$. Similarly, we have  $f^y\in G^1$ and $J(f)\in H^1$, hence, we obtain  the properties (a) and (b), and we may define a functional $I:L\to \mathbb R$ by
\[ I(f):=J(K(f)) \]
for all $f \in L$. It is immediate that $I$ is an integral and that the equation $I(g\otimes h)= J(g)K(h)$ is fulfilled. Thus we proved existence of a product integral.

Now let $I':L\to \mathbb R$ be  another integral fulfilling $I'(g\otimes h)=J(g)K(h)$. Then we have $I(f)=I'(f)$ for all $f \in G^1\otimes H^1$. For an arbitrary $f\in L$ we use the above $f_n\in G^1\otimes H^1$ to obtain
\begin{align*}
|I(f)-I'(f)| &\le I(|f-f_n|)+ |I(f_n)-I'(f_n)| + I'(|f_n-f|) \\
&\le n^{-1} I(g\otimes h)+ n^{-1}I'(g \otimes h).
\end{align*}
Taking the limit $n \to \infty$ yields $I(f)=I'(f)$, thus we achieve uniqueness of the product measure. In particular, setting $I'(f)= K(J(f))$ we obtain equality of the two iterated integrals, which finishes the proof.
\end{proof}

\noindent
For the proof of Theorem \ref{thm2} we resort to Stone's Fubini theorem. Stone \cite{Sto2} uses here the following strengthened version of the conditions (a) and (b):
\begin{enumerate}
\item[(a'')] $f^x \in H$, $f^y\in G$ for all $x\in X$, $y\in Y$,  
\item[(b'')] $K(f) \in G$, $J(f) \in H$.
\end{enumerate}
Now the Fubini-Stone theorem reads as follows.

\begin{prop}\label{prop2} Let $G\subset F(X)$, $H\subset F(Y)$ and $L\subset F(X\times Y)$ be Riesz spaces and   $J:G\to \mathbb R$, $K:H\to \mathbb R$ and  $I:L\to \mathbb R$ integrals. Then, if all $f\in L$ fulfil {\em (a'')}, {\em (b'')} and the equation $I(f)=J(K(f))$, then all $f \in L^1$ fulfil {\em (a')}, {\em (b')} and $I(f)=J(K(f))$.
\end{prop}

\noindent
For the proof see \cite{Sto2} or \cite[Theorem 7-2.I]{Tay}.

\begin{proof}[Proof of Theorem \ref{thm2}]
 Set $G=G^1$, $H=H^1$, $L=L(G^1\otimes H^1)$ and $I=J\otimes K$ in Proposition \ref{prop2}. Then, by Theorem \ref{thm1} the assumptions of the proposition are met. Thus, by its conclusion it follows $I(f)=J(K(f))$ for all $f\in L^1$. By symmetry, we also obtain $I(f)=K(J(f))$.
\end{proof}

\paragraph{Remark.} We note that with (a'') and (b'') instead of (a) and (b) the conclusion of the second part of Theorem~\ref{thm1} may fail, as seen from the following example:
Let $G=H$ be the Riesz space of all continuous piecewise linear functions on the interval $[0,1]$. Then the function $f(x,y):=y-x$ belongs to $G\otimes H$ and $f^+(x,y)=(y-x)^+$ is an element of $L(G\otimes H)$. Moreover, setting $J(g)= \int_0^1 g(x)\, dx$ for $g\in G$, we obtain an integral $J$ on $G$. Then $J(f^+)(y)=\int_0^1(y-x)^+ \, dx= y^2/2$, thus $J(f^+)$ does not belong to $H$.

\end{document}